\documentclass[twoside,12pt]{article}
\usepackage[all, 2cell, dvips]{xy}
\usepackage{latexsym,amsmath, amsfonts, graphics, epsf, epic}
\usepackage{amsthm}
\usepackage{amssymb}
\usepackage{amsopn}
\usepackage{amscd}
\usepackage{color}

\theoremstyle{plain}
\newtheorem{thm}{Theorem}[section]
\newtheorem{cor}[thm]{Corollary}
\newtheorem{lem}[thm]{Lemma}
\newtheorem{prop}[thm]{Proposition}

\theoremstyle{definition}

\newtheorem{remark}[thm]{Remark}
\newtheorem{example}[thm]{Example}
\newtheorem{defn}[thm]{Definition}

\usepackage{amssymb}

\def\lm{\lambda}

\def\al{\alpha}

\def\l.l.o.{\it l.l.o}

\def\chiup{\raise 2pt\hbox{$\chi$}}

\title{Asymptotics of Young tableaux in the $(k,\ell)$ hook}  

\author{A. Berele and A. Regev}

\begin{document}


\maketitle

\medskip

{\bf Abstract}. The asymptotics of the "$(k,\ell)$ hook" sums
$S_{k,\ell}^{(2z)}(n)$
(see~\eqref{the.sums.1}) were
calculated in~\cite{berele}.
It was recently realized that in~\cite[Section 7]{berele} there are few misprints and certain confusion with regard to the notations,  so that
the precise asymptotics of $S_{k,\ell}^{(2z)}(n)$
is not clear. Here we add more details and carefully repeat these
calculations, which lead to explicit values for the asymptotics of $S_{k,\ell}^{(2z)}(n).$


\section{Introduction}
Let $\lm$ be a partition and denote by $f^\lm$ the number of
standard Young tableaux (SYT) of shape $\lm$.
 Let $H(k,\ell;n)$ denote the partitions of $n$ in
the $(k,\ell)$ hook, namely
$H(k,\ell;n)=\{\lm=(\lm_1,\lm_2,\ldots)\mid \lm\vdash
n~\mbox{and}~ \lm_{k+1}\le \ell\}$.  In~\cite[Section 7]{berele}
we computed the asymptotics, as $n$ goes to infinity, of the sums
\begin{eqnarray}\label{the.sums.1}
S_{k,\ell}^{(2z)}(n)=\sum _{\lm\in
H(k,\ell;n)}\left(f^\lm\right)^{2z}.
\end{eqnarray}
That asymptotics has the form
\[
S_{k,\ell}^{(2z)}(n)\simeq a(k,\ell,2z)\cdot \left(\frac{1}{n}  \right)^{g(k,\ell,2z)}\cdot
(k+\ell)^{2zn}
\]
for some functions $a(k,\ell,2z)$ and $g(k,\ell,2z)$.

\medskip
It was recently realized that in~\cite[Section 7]{berele} there are few misprints and some confusion with the notations,  so that
the precise value of the constant term $a(k,\ell,2z)$ in that calculation is not clear.
Here we add more details and carefully repeat these
calculations, which lead to explicit values for the asymptotics of $S_{k,\ell}^{(2z)}(n), $
namely the explicit expression for the functions $a(k,\ell,2z)$ and $g(k,\ell,2z)$.
This is  Theorem~\ref{intro1} below (see also Theorem~\ref{final1}).

\medskip
Note that $\Gamma$ here is the gamma function.

\begin{thm}\label{intro1}
As $n$ goes to infinity,
\[
S_{k,\ell}^{(2z)}(n)\simeq\
a(k,\ell,2z)\cdot\left(\frac{1}{n}\right)^{g(k,\ell,2z)}\cdot
(k+\ell)^{2zn},
\]
where
\[g(k,\ell,2z)=
\frac{1}{2}\cdot( z\cdot[k(k+1)+\ell(\ell+1)-2]-(k+\ell-1))
\]
and
\begin{eqnarray*}
a(k,\ell,2z)=~~~~~~~~~~~~~~~~~~~~~~~~~~~~~~~~~~~~~~~~~~~~~~~~~~~~~~~~~~~~~~~~~~~~~~~~~~~~~~
~~~~~~~~~~~~~~~~~~~~~~~~~~~~~~~~~~~~~~~~~~~~~~~\\
=\left[\left(\frac{1}{\sqrt{2\pi}}\right)^{k+\ell-1}\cdot\left(\frac{1}{2}
\right)^{k\ell}\cdot(k+\ell)^{(k^2+\ell^2)/2}\right]^{2z}\cdot~~~~~~~~~~~~~~~~~~~~~~~~~~~~~~~~~~~~
~~~~~~~~~~~~~~~~~~~~~~~~~~~~~~~~~~~~~~~~~~~~~~\\
\cdot\left(\frac{1}{k+\ell}
\right)^{\frac{1}{2}\cdot[z\cdot(k(k-1)+\ell(\ell-1))+k+\ell]}\cdot
\frac{1}{k!\cdot\ell!}\cdot
\sqrt{\frac{z}{\pi}} \cdot\left(\sqrt{2\pi}
\right)^{k+\ell}\cdot(2z)^{-\frac{1}{2}\cdot
(z\cdot[k(k-1)+\ell(\ell-1)]+k+\ell)}\cdot~~~~~~~~~~~~~~~~~~~~~~~~~~~~~~\\
\cdot(\Gamma \left(1+z\right))^{-k-\ell}\cdot
\prod_{i=1}^k\Gamma\left(1+z i\right)\cdot
\prod_{j=1}^\ell\Gamma\left(1+z
j\right).~~~~~~~~~~~~~~~~~~~~~~~~~~~
\end{eqnarray*}
\end{thm}
In the case $2z=1$ we have (see Theorem~\ref{final2} below)

\begin{thm}\label{intro2}

\[
S_{k,\ell}^{(1)}(n)\simeq
a(k,\ell,1)\cdot\left(\frac{1}{n}\right)^{\frac{1}{4}\cdot(k(k-1)+\ell(\ell-1))}\cdot
(k+\ell)^{n}
\]
where
\begin{eqnarray*}
a(k,\ell,1)=~~~~~~~~~~~~~~~~~~~~~~~~~~~~~~~~~~~~~~~~~~~~~~~~~~~~~~~~~~~~~~~~~~~~~~~~~
~~~~~~~~~~~~~~~~~~~~~~~~~~~~~~~\\
=\left(\frac{1}{ 2}
\right)^{k\ell-k-\ell}\cdot\left(\frac{1}{\sqrt \pi}
\right)^{k+\ell}\cdot(k+\ell)^{\frac{1}{4}\cdot[k(k-1)+\ell(\ell-1)]}\cdot
~~~~~~~~~~~~~~~~~~~~~~~~~~~~~~~~~~~~~~~~~~\\
\frac{1}{k!\cdot\ell!}\cdot
\prod_{i=1}^k\Gamma\left(1+ i/2\right)\cdot
\prod_{j=1}^\ell\Gamma\left(1+j/2\right).~~~~~~~~~~~~~~~~~~~~~~~~
\end{eqnarray*}

\end{thm}

For the evaluation of special cases note that
\[
\Gamma\left(\frac{3}{2}\right)=\frac{\sqrt\pi}{2}\qquad\mbox{and}\qquad\Gamma(1+x)=x\Gamma(x).
\]

In several cases the sums $S_{k,\ell}^{(2z)}(n)$ can be evaluated
and given by simple formulas which yield the corresponding
asymptotics directly --
independent of Theorem~\ref{intro1}. These are the cases of
$S_{1,1}^{(1)}(n)$ and $S_{1,1}^{(2)}(n)$, see
Section~\ref{special.cases}, where we compare and verify that in
these cases the direct asymptotics does agree with the asymptotics deduced from
Theorem~\ref{intro1}.
It is also possible to use, say, "Mathematica"  to verify the validity of Theorem~\ref{intro1} in few special cases, see Section~\ref{math1}.

\section{Preliminaries}

\subsection{Recalling the "strip" case~\cite{regev}}

Theorem~\ref{intro1} is a hook generalization of the following
"strip" theorem~\ref{intro3},
see~\cite[Corollary 4.4]{regev}. We remark that, even though
Theorem~\ref{intro1} is proved under the assumption that both
$k,\ell\ge 1$, nevertheless it reduces to Theorem~\ref{intro3}
when one substitutes $\ell=0$.
\begin{thm}\label{intro3}
Let
\[
S^{(2z)}_k(n)=\sum_{\lm\in H(k,0;n)}(f^\lm)^{2z},
\]
then, as $n$ goes to infinity,
\[
S^{(2z)}_k(n)\simeq
a(k,2z)\cdot\left(\frac{1}{n}\right)^{g(k,z)}\cdot k^{2z n}
\]
where
\[
g(k,z)=\frac{1}{2}\cdot (z(k^2+k-2)-(k-1))
\]
and
\begin{eqnarray*}
a(k,2z)= \left[\left(\frac{1}{\sqrt{2\pi}} \right)^{k-1}\cdot
k^{k^2/2}\right]^{2z}
\cdot\left(\frac{1}{k}\right)^{[zk(k-1)+k]/2}\cdot
\frac{1}{k!}\cdot\sqrt{\frac{z}{\pi}}\cdot (2\pi)^{k/2}\cdot
(2z)^{-(zk(k-1)+k)/2}
\cdot~~~~~~~~~~~~~~~~~~~~~~~~~~~~~\\
\cdot \Gamma(1+z)^{-k}\cdot\prod_{j=1}^k\Gamma(1+z
j).~~~~~~~~~~~~~~~~~~~~~~~~~
\end{eqnarray*}
\end{thm}
\begin{remark}
One of the main tools in proving Theorem~\ref{intro3}
in~\cite{regev} was the computation of the asymptotics of a single
$f^\nu$ where $\nu\in H(k,0;n)$. Let
$\nu=(\nu_1,\ldots,\nu_k)\vdash n$, write
$\nu_i=\frac{n}{k}+a_i\sqrt n$. By~\cite[(F.1.1)]{regev}, when the
$a_i$s are bounded
we have:
\begin{eqnarray}\label{nu0.1.2}
f^\nu\simeq \gamma_k\cdot D_k(a_1,\ldots,a_k)\cdot e^{-(k/2)(\sum
a_i^2)}\cdot \left(\frac{1}{n}\right )^{(k-1)(k+2)/4}\cdot k^{n},
\end{eqnarray} where
\begin{eqnarray}\label{nu022}
\gamma_k=\left(\frac{1}{\sqrt{2\pi}}\right)^{k-1}\cdot k^{k^2/2}
\qquad\mbox{and}\qquad D_k(a_1,\ldots,a_k)=\prod_{1\le i<j\le
k}(a_i-a_j).
\end{eqnarray}
We apply~\eqref{nu0.1.2}
in what follows.
\end{remark}

\subsection{Preliminaries for the hook case}\label{section.2.1}

We assume that both $k,\ell\ge 1$ and we follow the notations
of~\cite{berele} from  Section 7.9 on. We assume that $\lm\in
H(k,\ell;n)$ and $\lm_k\ge \ell$, namely $\lm$ contains the
 $k\times \ell$ rectangle $R_{k,\ell}$. Thus $\lm$ is
made of the partitions $\nu,\mu'$ and of the rectangle
$R_{k,\ell}$. We have $\lm\vdash n$, $\nu\vdash n_k$, $\mu\vdash
n_\ell$, $R_{k,\ell}\vdash k\ell$, $n=n_k+n_\ell+k\ell$, $\bar
n=n-k\ell$.

\medskip
By assumption $n_k\simeq n\cdot \frac{k}{k+\ell}$ and
$n_\ell\simeq n\cdot \frac{\ell}{k+\ell}$. Now $\nu=(\nu_1,\ldots
\nu_k)$, $\mu=(\mu_1,\ldots,\mu_\ell)$, and
\begin{eqnarray}\label{star1}
\nu_i=\frac{n_k}{k}+a_i\sqrt{n_k},\qquad
\mu_j=\frac{n_\ell}{\ell}+b_j\sqrt{n_\ell}~.
\end{eqnarray}
Also
\begin{eqnarray}\label{star2}
\nu_i=\frac{\bar n}{k+\ell}+\al_i\sqrt{\bar n}, \qquad
\mu_j=\frac{\bar n}{k+\ell}+\beta_j\sqrt{\bar n},
\end{eqnarray}
and we write $\al_1+\cdots+\al_k=\al$ and
$\beta_1+\cdots+\beta_\ell=\beta$. It follows that $\beta=-\al$.
The transition from the $a_i,b_j$ to the $\al_i, \beta_j $ is
given by
\[
n_k=\frac{k\bar n}{k+\ell}+\al\sqrt{\bar n}, \qquad\qquad
n_\ell=\frac{\ell\bar n}{k+\ell}-\al\sqrt{\bar n},
\qquad\mbox{hence, since $\bar n\to\infty$,}
\]
\begin{eqnarray}\label{51}
a_i =\left(\al_i-\frac{\al}{k}
\right)\cdot\left(\frac{k}{k+\ell}+\frac{\al}{\sqrt {\bar n}}
\right)^{-1/2}\simeq \left(\al_i-\frac{\al}{k}
\right)\cdot\left(\frac{k+\ell}{k} \right)^{\frac{1}{2}}
\end{eqnarray}
and the difference $l.h.s.-r.h.s$ tends to zero as $n\to\infty$. Similarly
\begin{eqnarray}\label{52}
 b_j =\left(\beta_j+\frac{\al}{\ell}
\right)\cdot\left(\frac{\ell}{k+\ell}-\frac{\al}{\sqrt {\bar n}}
\right)^{-1/2}\simeq\left(\beta_j+\frac{\al}{\ell}
\right)\cdot\left(\frac{k+\ell}{\ell} \right)^{\frac{1}{2}}
\end{eqnarray}
and $~l.h.s.-r.h.s\longrightarrow 0$.

\section{Asymptotic of a single $f^\lm$}

We refer now to~\cite[Section 7]{berele}. Up to Lemma 7.15 there,
including that lemma, every detail was checked and verified, and we proceed
from that point.

\medskip
The hook formula yields the factorisation $f^\lm=A_1\cdot A_2\cdot
A_3\cdot A_4 $ (see~\cite[7.14]{berele}) where

\medskip
$A_1=n!/\bar n!\simeq n^{k\ell},  $

\medskip
$A_2=1/(\prod _R h_{ij})\simeq ((k+\ell)/(2n))^{k\ell}$,

\medskip
$A_3=f^\nu\cdot f^\mu$, ~~~ and
\[
A_4=\frac{\bar n!}{n_k!\cdot n_\ell !}\simeq
\frac{1}{\sqrt{2\pi}}\cdot (k+\ell)^{-k\ell}\cdot
\frac{k+\ell}{\sqrt{k\ell}}\cdot \frac{1}{\sqrt{n}}\cdot
\frac{(k+\ell)^n}{k^{n_k}\cdot \ell^{n_\ell}}\cdot
e^{-\al^2(k+\ell)^2/(2k\ell)}.~~~~~~~~~~~~~~~~~~~~~~~~~~~~
\]
Note that $A_1\cdot A_2\simeq ((k+\ell)/2)^{kl}$, so
\begin{eqnarray}\label{nu0.1}
f^\lm=A_1\cdot A_2\cdot A_3\cdot
A_4\simeq\frac{1}{\sqrt{2\pi}}\cdot \frac{1}{2^{k\ell}} \cdot
\frac{k+\ell}{\sqrt{k\ell}}\cdot \frac{1}{\sqrt{n}}\cdot
\frac{(k+\ell)^n}{k^{n_k}\cdot \ell^{n_\ell}}\cdot
e^{-\al^2(k+\ell)^2/(2k\ell)}\cdot f^\nu\cdot f^\mu.
\end{eqnarray}
We analyze $f^\nu$ and  $f^\mu$. By~\eqref{nu0.1.2} and by
of~\cite[(F.1.1)]{regev},
\begin{eqnarray}\label{nu1}
f^\nu\simeq \gamma_k\cdot D_k(a_1,\ldots,a_k)\cdot e^{-(k/2)(\sum
a_i^2)}\cdot \left(\frac{1}{n_k}\right )^{(k-1)(k+2)/4}\cdot
k^{n_k}, \end{eqnarray} where
$\gamma_k=({1}/{\sqrt{2\pi}})^{k-1}\cdot k^{k^2/2}.$
We make the transition from the $a_i$ to the $\al_i$.
By~\eqref{51} $a_i-a_j\simeq(\al_i-\al_j)\cdot\sqrt{(k+\ell)/k}$, hence
\begin{eqnarray}\label{nu3}
D_k(a_1,\ldots,a_k)\simeq\left(\frac{k+\ell}{k}\right)^{\frac{k(k-1)}{4}}\cdot
D_k(\al_1,\ldots,\al_k),
\end{eqnarray}
and similarly by~\eqref{52}
\begin{eqnarray}\label{nu033}
D_\ell(b_1,\ldots,b_\ell)\simeq\left(\frac{k+\ell}{\ell}\right)^{\frac{\ell(\ell-1)}{4}}\cdot
D_\ell(\beta_1,\ldots,\beta_\ell).
\end{eqnarray}
{\bf Claim:}
\begin{eqnarray}\label{nu44}
e^{-\frac{k}{2}\sum a_i^2}\simeq
e^{-\frac{(k+\ell)}{2}\sum\al_i^2}\cdot
e^{\left(\frac{k+\ell}{2k}\right)\cdot\al^2}
\end{eqnarray}
\begin{proof}
By~\eqref{51} \[a_i^2\simeq \frac{k+\ell}{k}\cdot
\left(\al_i^2-\frac{2\al}{k}\al_i+\frac{\al^2}{k^2}\right)\qquad\mbox{and
$~l.h.s-r.h.s\longrightarrow 0$  as $n\to\infty$.}\]  Since $\sum \al_i=\al$, we
have
\[
\sum a_i^2\simeq \frac{k+\ell}{k}\cdot
\sum\al_i^2-\frac{k+\ell}{k^2}\cdot\al^2\quad\mbox{and
$~l.h.s-r.h.s\longrightarrow 0$, \qquad so}
\]
\[
\quad\frac{k}{2}\cdot \sum a_i^2\simeq \frac{k+\ell}{2}\cdot
\sum\al_i^2-\frac{k+\ell}{2k}\cdot\al^2\quad\mbox{and
$~l.h.s-r.h.s\longrightarrow 0$  as $n\to\infty$,}
\]
and this implies~\eqref{nu44}.
\end{proof}

\begin{cor}\label{nu5}
By~\eqref{nu1},~\eqref{nu3} and~\eqref{nu44}
\begin{eqnarray}\label{nu5}
f^\nu\simeq\gamma_k\cdot
\left(\frac{k+\ell}{k}\right)^{k(k-1)/4}\cdot
D_k(\al_1,\ldots,\al_k)\cdot e^{-\frac{k+\ell}{2}\sum \al_i^2}
\cdot e^{+\frac{k+\ell}{2k}\al^2} \cdot\left(\frac{1}{n_k}
\right)^{(k-1)(k+2)/4}\cdot k^{n_k},
\end{eqnarray}
and similarly (recall that $\beta=-\al$)
\begin{eqnarray}\label{nu6}
f^\mu\simeq\gamma_\ell \cdot \left(\frac{k +\ell}{\ell
}\right)^{\ell (\ell -1)/4}\cdot D_\ell (\beta_1,\ldots,\beta_\ell
)\cdot e^{-\frac{k+\ell }{2}\sum \beta_j^2} \cdot e^{+\frac{k
+\ell}{2\ell }\al^2} \cdot\left(\frac{1}{n_\ell } \right)^{(\ell
-1)(\ell +2)/4}\cdot \ell ^{n_\ell }.
\end{eqnarray}
\end{cor}

Plugging~\eqref{nu5} and~\eqref{nu6} into~\eqref{nu0.1} gives
\begin{cor}\label{late1}
\begin{eqnarray*}\label{nu4}
f^\lm\simeq \frac{1}{\sqrt{2\pi}}\cdot\left(\frac{1}{2}
\right)^{k\ell}\cdot\frac{k+\ell}{\sqrt{k\ell}}\cdot\frac{1}{\sqrt{n}}\cdot
e^{-\frac{(k+\ell)^2}{2k\ell}\cdot\al^2}\cdot
\frac{(k+\ell)^n}{k^{n_k}\ell^{n_\ell}}\cdot f^\nu\cdot
f^\mu\simeq~~~~~~~~~~~~~~~~~~~~~~~~~~~~~~~~~~~~~~\\
\simeq\frac{1}{\sqrt{2\pi}}\cdot\left(\frac{1}{2}
\right)^{k\ell}\cdot\frac{k+\ell}{\sqrt{k\ell}}\cdot\frac{1}{\sqrt{n}}\cdot
e^{-\frac{(k+\ell)^2}{2k\ell}\cdot\al^2}\cdot
\frac{(k+\ell)^n}{k^{n_k}\ell^{n_\ell}}\cdot~~~~~~~~~~~~~~~~~~~~~~~~~~~~~~~~~~~~~~~~~~~~~~~~~~~\\
\gamma_k\cdot \left(\frac{k+\ell}{k}\right)^{k(k-1)/4}\cdot
D_k(\al_1,\ldots,\al_k)\cdot e^{-\frac{k+\ell}{2}\sum \al_i^2}
\cdot e^{+\frac{k+\ell}{2k}\al^2} \cdot\left(\frac{1}{n_k}
\right)^{(k-1)(k+2)/4}\cdot k^{n_k}~~~~~~~\\
\gamma_\ell \cdot \left(\frac{k +\ell}{\ell }\right)^{\ell (\ell
-1)/4}\cdot D_\ell (\beta_1,\ldots,\beta_\ell )\cdot
e^{-\frac{k+\ell }{2}\sum \beta_j^2} \cdot e^{+\frac{k
+\ell}{2\ell }\al^2} \cdot\left(\frac{1}{n_\ell } \right)^{(\ell
-1)(\ell +2)/4}\cdot \ell ^{n_\ell }~~~~~~~
\end{eqnarray*}
\end{cor}

Collecting terms in Corollary~\ref{late1}, we proved the following
theorem (which is~\cite[ Theorem
7.16]{berele}).

\begin{thm}\label{nu32}
With the notations of Section~\ref{section.2.1} we have
\begin{eqnarray*}
f^\lm\simeq c(k,\ell)\cdot D_k(\al_1,\ldots,\al_k)\cdot
D_\ell(\beta_1,\ldots,\beta_\ell)\cdot
e^{-\frac{k+\ell}{2}(\sum\al_i^2+\sum\beta_j^2)} \cdot\left(
\frac{1}{n}\right)^{\theta(k,\ell)}\cdot(k+\ell)^n
\end{eqnarray*}
with
\[
c(k,\ell)=\left(\frac{1}{\sqrt{2\pi}}\right)^{k+\ell-1}\cdot\left(\frac{1}{2}
\right)^{k\ell}\cdot(k+\ell)^{(k^2+\ell^2)/2} \] and ~~
\[\theta
(k,\ell)=\frac{1}{4}\cdot[k(k+1)+\ell(\ell+1)-2].
\]
\end{thm}

\section{Asymptotics for the sums $S_{k,\ell}^{(2z)}(n)$}
By~\eqref{the.sums.1}
\[
S_{k,\ell}^{(2z)}(n)=\sum_{\lm\in H(k,\ell;n)}(f^\lm)^{2z}.
\]
As in~\cite[Theorem 7.18]{berele}  (but with the additional factor
$e^{-\frac{(k+\ell)^2}{k\ell} \cdot u^2}$ in the integral), deduce

\begin{thm}\label{sum1}
With the notations of Theorem~\ref{nu32}, as $n$ goes to infinity
we have
\[
S_{k,\ell}^{(2z)}(n)\simeq\left[c(k,\ell)\cdot\left(\frac{1}{n}
\right)^{\theta (k,\ell)}\cdot (k+\ell)^n   \right]^{2z}\cdot
(\sqrt n)^{k+\ell-1}\cdot I(k,\ell,2z),
\]
where $\sum x_i=u$, where
\begin{eqnarray}\label{i}
 I(k,\ell,2z)=\int_{P(k,\ell)}\left[D_k(x)\cdot D_\ell(y)\cdot
e^{-\frac{k+\ell}{2}(\sum x_i^2+\sum y_j^2)}\right]^{2z}
\;d^{(k+\ell-1)}(x,y),
\end{eqnarray}
and where $P(k,\ell)\subset \mathbb{R}^{k+\ell}$ is the domain
\[
P(k,\ell)=\{(x_1,\ldots,x_k;y_1,\dots ,y_\ell)\mid x_1\ge\cdots\ge
x_k;\;y_1\ge\ldots\ge y_\ell, ~\sum x_i+\sum y_j=0 \}.
\]
\end{thm}
Note that $\sum y_j=-u$ since $\sum x_i=u$ and $\sum x_i+\sum y_j=0$.

\subsection{The evaluation of $I(k,\ell,2z)$}
Let
\[
A_{k,\ell}^{(2z)}(u)=\int_{\Omega(k,u) }\left[D_k(x)\cdot
e^{-\frac{k+\ell}{2}\cdot\sum
x_i^2} \right ]^{2z}\;d^{(k-1)}(x)
\]
where $\Omega(k,u)\subset\mathbb{R}^k$, $~~\Omega(k,u)
=\{x_1\ge\cdots\ge x_k\mid \sum x_i=u\}$.

\medskip
Similarly let
\[
B_{k,\ell}^{(2z)}(-u)=\int_{\Omega(\ell,-u) }\left[D_\ell(y)\cdot
e^{-\frac{k+\ell}{2}\cdot\sum y_j^2} \right
]^{2z}\;d^{(\ell-1)}(y)
\]
where $\Omega(\ell,-u)\subset\mathbb{R}^\ell$,
~~$\Omega(\ell,-u)=\{y_1\ge\cdots\ge y_\ell\mid \sum y_j=-u\}$.

\medskip
We clearly have
\begin{lem}
The integral~\eqref{i} satisfies
\begin{eqnarray}\label{int1}
I(k,\ell,2z)=\int_{-\infty}^{\infty} A_{k,\ell}^{(2z)}(u)\cdot
B_{k,\ell}^{(2z)}(-u)\;du.
\end{eqnarray}
\end{lem}

\subsubsection{Evaluating $A_{k,\ell}^{(2z)}(u)$ and $B_{k,\ell}^{(2z)}(-u)$}

Recall that $\sum x_i=u$ and make the substitution
$x_i'=x_i-\frac{u}{k}$ (similarly $y_j'=y_j+\frac{u}{\ell}$), then
$\sum x_i'=\sum y_j'=0$; also $~D_k(x')=D_k(x)$ and
$D_\ell(y')=D_\ell(y)$. Also the Jacobians equal 1.

\medskip
Now $\sum x_i^2=\frac{u^2}{k}+\sum {x_i'}^2$ hence
\[
e^{-\frac{k+\ell}{2}\cdot\sum
x_i^2}=e^{-\frac{k+\ell}{2}\cdot\sum {x_i'}^2}\cdot
e^{-\frac{k+\ell}{2k}\cdot u^2}.
\]
It follows that \[
A_{k,\ell}^{(2z)}(u)=e^{-\frac{(k+\ell)z}{k}\cdot u^2} \cdot
I^{(2z)}_k
\]
where
\begin{eqnarray}\label{ii}
I^{(2z)}_k=\int_{\Omega'(k) }\left[D_k(x')\cdot
e^{-\frac{k+\ell}{2}\cdot\sum {x_i'}^2} \right
]^{2z}\;d^{(k-1)}(x')
\end{eqnarray}
and where $\Omega'(k)=\{x'_1\ge\cdots\ge x'_k\mid \sum x'_i=0\}$.

\medskip
Similarly
\[ B_{k,\ell}^{(2z)}(-u)=e^{-\frac{(k+\ell)z}{\ell}\cdot u^2} \cdot
I^{(2z)}_\ell
\]
where
\[
I^{(2z)}_\ell=\int_{\Omega'(\ell) }\left[D_\ell(y')\cdot
e^{-\frac{k+\ell}{2}\cdot\sum {y_j'}^2} \right
]^{2z}\;d^{(\ell-1)}(y')
\]
and where $\Omega'(\ell)=\{y'_1\ge\cdots\ge y'_\ell\mid \sum
y'_j=0\}$.

\medskip
By~\eqref{int1} we have
\[
I(k,\ell,2z)=I^{(2z)}_k\cdot
I^{(2z)}_\ell\cdot\int_{-\infty}^{\infty}
e^{-\frac{(k+\ell)z}{k}\cdot u^2}\cdot
e^{-\frac{(k+\ell)z}{\ell}\cdot u^2}\;du .
\]
Note that
\[
J:=\int_{-\infty}^{\infty} e^{-\frac{(k+\ell)z}{k}\cdot u^2}\cdot
e^{-\frac{(k+\ell)z}{\ell}\cdot u^2}\;du =\int_{-\infty}^{\infty}
e^{-\frac{z(k+\ell)^2}{k\ell}\cdot u^2}\;du.
\]
Thus we proved
\begin{lem}\label{int2}
\[
I(k,\ell,2z)=\int_{-\infty}^{\infty}A^{(2z)}(u)\cdot
B^{(2z)}(-u)\;du = I^{(2z)}_k \cdot I^{(2z)}_\ell
\cdot\int_{-\infty}^{\infty} e^{-\frac{z(k+\ell)^2}{k\ell}\cdot
u^2} du.
\]
\end{lem}

{\bf Calculate:}

Let $J=\int_{-\infty}^{\infty} e^{-\frac{z(k+\ell)^2}{k\ell}\cdot
u^2} du$ and denote $r=\sqrt{\frac{z(k+\ell)^2}{k\ell}}$. Since
$\int_{-\infty}^{\infty}e^{-r^2\cdot u^2}du=
\sqrt\pi/r,$ hence
\[
J=\sqrt{\frac{k\cdot\ell\cdot{\pi}}{z(k+\ell)^2}}.
\]
This implies

\begin{cor}\label{int9}
\begin{eqnarray*}
I(k,\ell,2z)=\sqrt{\frac{k\cdot\ell\cdot{\pi}}{z(k+\ell)^2}}\cdot
I^{(2z)}_k \cdot
I^{(2z)}_\ell=~~~~~~~~~~~~~~~~~~~~~~~~~~~~~~~~~~~~~~~~~~~~~~~~~~~~~~~~~~~~~~~~~~~~~~~~~~
\\=\sqrt{\frac{k\cdot\ell\cdot{\pi}}{z(k+\ell)^2}}\cdot\int_{\Omega'(k)
}\left[D_k(x')\cdot e^{-\frac{k+\ell}{2}\cdot\sum {x_i'}^2} \right
]^{2z}\;d^{(k-1)}(x')\cdot \int_{\Omega'(\ell)
}\left[D_\ell(y')\cdot e^{-\frac{k+\ell}{2}\cdot\sum {y_j'}^2}
\right
]^{2z}\;d^{(\ell-1)}(y')\\=\sqrt{\frac{k\cdot\ell\cdot{\pi}}{z(k+\ell)^2}}\cdot\int_{\Omega'(k)
}\left[D_k(x)\cdot e^{-\frac{k+\ell}{2}\cdot\sum {x_i}^2} \right
]^{2z}\;d^{(k-1)}(x)\cdot \int_{\Omega'(\ell)
}\left[D_\ell(y)\cdot e^{-\frac{k+\ell}{2}\cdot\sum {y_j}^2}
\right ]^{2z}\;d^{(\ell-1)}(y)
\end{eqnarray*}
(in the last term we replaced $x'$ by $x$ and $y'$ by $y$).
\end{cor}
\subsection{Evaluating $I^{(2z)}_k$ and $I^{(2z)}_\ell$}

The key to the following evaluation is the celebrated Selberg
integral~\cite{forrester},~\cite{selberg}. Let
\begin{eqnarray}\label{int4}
I(s,\beta):=\int_{\Omega '(s)}\left[|D_k(x)|\cdot
e^{-\frac{1}{2}\sum x_i^2}\right]^\beta \;dx_1\cdots dx_{s-1}
\end{eqnarray}
and let
\begin{eqnarray}\label{int04}
\Psi_s^{(\beta)}=(\sqrt{2\pi})^s\cdot \beta^{-s/2-\beta
s(s-1)/4}\cdot \left[\Gamma
\left(1+\frac{1}{2}\beta\right)\right]^{-s}\cdot
\prod_{j=1}^s\Gamma\left(1+\frac{1}{2}\beta j\right),
\end{eqnarray}
then it follows from the Selberg integral that
\begin{eqnarray}\label{int3}
I(s,\beta)=\frac{1}{s!}\cdot\frac{1}{\sqrt s}\cdot
\sqrt{\frac{\beta}{2\pi}}\cdot \Psi_s^{(\beta)},
\end{eqnarray}
see~\cite[(F.4.1) and (F.4.3)]{regev} . For elementary proofs of
both the Mehta and Selberg integrals see~\cite{garsia}.

\medskip
Thus, with $k=s$ and $\beta=2z$ we deduce
\begin{lem}\label{int6}
\begin{eqnarray*}
I(k,2z)=\frac{1}{k!}\cdot\frac{1}{\sqrt k}\cdot
\sqrt{\frac{2z}{2\pi}}\cdot \Psi_k^{(2z)}=~~~~~~~~~~~~~~~~~~~~~~~~~~~~~~~~~~~~~~~~~~~~~~~~~~~~~~~~~~~~~~~~~~~~~~~~~~\\
~~~~~~~~~~~=\frac{1}{k!}\cdot\frac{1}{\sqrt k}\cdot
\sqrt{\frac{z}{\pi}}\cdot (\sqrt{2\pi})^k\cdot (2z)^{-k/2-2z
k(k-1)/4}\cdot \left[\Gamma \left(1+z\right)\right]^{-k}\cdot
\prod_{j=1}^k\Gamma\left(1+z j\right).~~~~~~~~~~~~~~~
\end{eqnarray*}
\end{lem}

We rewrite~\eqref{ii} as
\begin{eqnarray}\label{iii}
I^{(2z)}_k=\int_{\Omega'(k) }\left[D_k(x_1,\ldots,x_k)\cdot
e^{-\frac{k+\ell}{2}\cdot\sum {x_i}^2} \right
]^{2z}\;d^{(k-1)}(x),
\end{eqnarray}
and make the transition from ~\eqref{int4} to~\eqref{iii} as follows.
\begin{lem}\label{int7}
With $I(k,2z)$ given by~\eqref{int4},
\[
I^{(2z)}_k=\left(\frac{1}{{k+\ell}}\right)^{\frac{k(k-1)z+k-1}{2}}\cdot
I(k,2z)
\]
\end{lem}
\begin{proof}
In~\eqref{int4}
substitute $x_i=\sqrt{\frac{1}{k+\ell}}\cdot
y_i$, then
\[
D_k(x)=\left(\frac{1}{\sqrt{k+\ell}}\right)^{\frac{k(k-1)}{2}}D_k(y)=
\left(\frac{1}{{k+\ell}}\right)^{\frac{k(k-1)}{4}}D_k(y)\] and
\[d^{(k-1)}(x)=\left(\frac{1}{\sqrt{k+\ell}}\right)^{k-1}d^{(k-1)}(y).
\]
Thus
\begin{eqnarray*}
I^{(2z)}_k=\left(\frac{1}{{k+\ell}}\right)^{\frac{k(k-1)z+k-1}{2}}\cdot
\int_{\Omega'(k)}\left[D_k(y_1,\ldots,y_k)\cdot
e^{-\frac{1}{2}\cdot\sum
y_i^2}\right]^{2z}d^{(k-1)}(y)=~~~~~~~~~~~~~~~~~~~~~~~~~~~~~~~~~~~\\
=\left(\frac{1}{{k+\ell}}\right)^{\frac{k(k-1)z+k-1}{2}}\cdot
I(k,2z)~~~~~~~~~~~~~~~~~~~~~~
\end{eqnarray*}
\end{proof}
Together with Lemma~\ref{int6} it implies
\begin{cor}\label{int8}
\[
I^{(2z)}_k=\left(\frac{1}{{k+\ell}}\right)^{\frac{k(k-1)z+k-1}{2}}\cdot\frac{1}{k!}\cdot\frac{1}{\sqrt
k}\cdot \sqrt{\frac{z}{\pi}}\cdot (\sqrt{2\pi})^k\cdot
(2z)^{-k/2-z k(k-1)/2}\cdot \left[\Gamma
\left(1+z\right)\right]^{-k}\cdot \prod_{j=1}^k\Gamma\left(1+z
j\right).
\]
Similarly
\[
I^{(2z)}_\ell=\left(\frac{1}{{k+\ell}}\right)^{\frac{\ell(\ell-1)z+\ell-1}{2}}\cdot\frac{1}{\ell!}\cdot\frac{1}{\sqrt
\ell}\cdot \sqrt{\frac{z}{\pi}}\cdot (\sqrt{2\pi})^\ell\cdot
(2z)^{-\ell/2-z \ell(\ell-1)/2}\cdot \left[\Gamma
\left(1+z\right)\right]^{-\ell}\cdot
\prod_{j=1}^\ell\Gamma\left(1+z j\right).
\]
\end{cor}

\medskip
By Corollary ~\ref{int9} we get that
\begin{eqnarray*}
I^{(2z)}_k \cdot I^{(2z)}_\ell=
~~~~~~~~~~~~~~~~~~~~~~~~~~~~~~~~~~~~~~~~~~~~~~~~~~~~~~~~~~~~~~~~~~~~~~~~~~~~~~~~~~~~~~~~~~~~~~~~~~~~~~~~~~~~~~~~~~\\
=\left(\frac{1}{{k+\ell}}\right)^{\frac{k(k-1)z+k-1}{2}}\cdot\frac{1}{k!}\cdot\frac{1}{\sqrt
k}\cdot \sqrt{\frac{z}{\pi}}\cdot (\sqrt{2\pi})^k\cdot
(2z)^{-k/2-z k(k-1)/2}\cdot \left[\Gamma
\left(1+z\right)\right]^{-k}\cdot \prod_{i=1}^k\Gamma\left(1+z
i\right)\cdot~~~~~\\
\left(\frac{1}{{k+\ell}}\right)^{\frac{\ell(\ell-1)z+\ell-1}{2}}\cdot\frac{1}{\ell!}\cdot\frac{1}{\sqrt
\ell}\cdot \sqrt{\frac{z}{\pi}}\cdot (\sqrt{2\pi})^\ell\cdot
(2z)^{-\ell/2-2z \ell(\ell-1)/2}\cdot \left[\Gamma
\left(1+z\right)\right]^{-\ell}\cdot
\prod_{j=1}^\ell\Gamma\left(1+z j\right)=~~~~~\\
= \left(\frac{1}{k+\ell}
\right)^{\frac{1}{2}\cdot[z\cdot(k(k-1)+\ell(\ell-1))+k+\ell-2]}\cdot
\frac{1}{k!\cdot\ell!}\cdot\frac{1}{\sqrt{k\ell}}\cdot\frac{z}{\pi}\cdot
\left(\sqrt{2\pi} \right)^{k+\ell}\cdot(2z)^{-\frac{1}{2}\cdot
(z\cdot[k(k-1)+\ell(\ell-1)]+k+\ell)}\cdot~~~~~~~~~~~~~\\
\cdot (\Gamma \left(1+z\right))^{-k-\ell}\cdot
\prod_{i=1}^k\Gamma\left(1+z i\right)\cdot
\prod_{j=1}^\ell\Gamma\left(1+z j\right).~~~~~~
\end{eqnarray*}
Therefore
\begin{eqnarray*}
I(k,\ell,2z) =~~~~~~~~~~~~~~~~~~~~~~~~~~~~~~~~~~~~~~~~~~~~~~~~~~~~~~~~~~~~~~~~~~~~~~~~~~~~~~~~~~
~~~~~~~~~~~~~~~~~~~~~~~~~~~~~~~~~~~~~~~~~~\\
=\frac{1}{k+\ell}\cdot\sqrt{\frac{k\cdot\ell\cdot{\pi}}{z}}\cdot~~~~~~~~~~~~~~~~~~~~~~~~~~~~~~~~~~~~~
~~~~~~~~~~~~~~~~~~~~~~~~~~~~~~~~~~~~~~~~~~~~~~~~~~~~~~~~~~~~~~~~~~~~~~~~~~~~\\
\cdot \left(\frac{1}{k+\ell}
\right)^{\frac{1}{2}\cdot[z\cdot(k(k-1)+\ell(\ell-1))+k+\ell-2]}\cdot
\frac{1}{k!\cdot\ell!}\cdot
\frac{1}{\sqrt{k\ell}}\cdot\frac{z}{\pi}\cdot \left(\sqrt{2\pi}
\right)^{k+\ell}\cdot(2z)^{-\frac{1}{2}\cdot
(z\cdot[k(k-1)+\ell(\ell-1)]+k+\ell)}\cdot~~~~~~~~~~~~~~~~~~\\
\cdot (\Gamma \left(1+z\right))^{-k-\ell}\cdot
\prod_{i=1}^k\Gamma\left(1+z i\right)\cdot
\prod_{j=1}^\ell\Gamma\left(1+z j\right)~~~~~~~~~~~~~~~~
\end{eqnarray*}
so, after cancellations,
\begin{eqnarray*}
I(k,\ell,2z) =\left(\frac{1}{k+\ell}
\right)^{\frac{1}{2}\cdot[z\cdot(k(k-1)+\ell(\ell-1))+k+\ell]}\cdot
\frac{1}{k!\cdot\ell!}\cdot
\sqrt{\frac{z}{\pi}} \cdot\left(\sqrt{2\pi}
\right)^{k+\ell}\cdot(2z)^{-\frac{1}{2}\cdot
(z\cdot[k(k-1)+\ell(\ell-1)]+k+\ell)}\cdot~~~~~~~~~~~~~~~~~~~~~~~~~~~~~~\\
(\Gamma \left(1+z\right))^{-k-\ell}\cdot
\prod_{i=1}^k\Gamma\left(1+z i\right)\cdot
\prod_{j=1}^\ell\Gamma\left(1+z
j\right).~~~~~~~~~~~~~~~~~~~~~~~~~~~~~
\end{eqnarray*}
Combined with Theorem~\ref{sum1} we have proved

\begin{thm}\label{final1}(This is also Theorem~\ref{intro1})
\[
S_{k,\ell}^{(2z)}\simeq\
a(k,\ell,2z)\cdot\left(\frac{1}{n}\right)^{g(k,\ell,2z)}\cdot
(k+\ell)^{2zn},
\]
where
\[g(k,\ell,2z)=
\frac{1}{2}\cdot( z\cdot[k(k+1)+\ell(\ell+1)-2]-(k+\ell-1))
\]
and
\begin{eqnarray*}
a(k,\ell,2z)=[c(k,\ell)]^{2z}\cdot
I(k,\ell,2z)=~~~~~~~~~~~~~~~~~~~~~~~~~~~~~~~~~~~~~~~~~~~~~~~~~~~~~~~~~~~~~~~~~~~~~~~~~~~~~~~~~~~~~~~~~~~~~~~~~~~~~~\\
=\left[\left(\frac{1}{\sqrt{2\pi}}\right)^{k+\ell-1}\cdot\left(\frac{1}{2}
\right)^{k\ell}\cdot(k+\ell)^{(k^2+\ell^2)/2}\right]^{2z}\cdot~~~~~~~~~~~~~~~~~~~~~~~~~~~~~~~~~~~~~~~~~~
~~~~~~~~~~~~~~~~~~~~~~~~~~~~~~~~~~~~~~~~~~~~~~\\
\cdot\left(\frac{1}{k+\ell}
\right)^{\frac{1}{2}\cdot[z\cdot(k(k-1)+\ell(\ell-1))+k+\ell]}\cdot
\frac{1}{k!\cdot\ell!}\cdot
\sqrt{\frac{z}{\pi}} \cdot\left(\sqrt{2\pi}
\right)^{k+\ell}\cdot(2z)^{-\frac{1}{2}\cdot
(z\cdot[k(k-1)+\ell(\ell-1)]+k+\ell)}\cdot~~~~~~~~~~~~~~~~~~~~~~~~~~~~~~\\
\cdot(\Gamma \left(1+z\right))^{-k-\ell}\cdot
\prod_{i=1}^k\Gamma\left(1+z i\right)\cdot
\prod_{j=1}^\ell\Gamma\left(1+z
j\right).~~~~~~~~~~~~~~~~~~~~~~~~~~~
\end{eqnarray*}
\end{thm}
\section{Some special cases}\label{special.cases}

\subsection{The case $2z=1$}
Here
\[
g(k,\ell,1)=\frac{1}{4}\cdot(k(k-1)+\ell(\ell-1)).
\]
We calculate $a(k,\ell,1)$.  Recall that $\Gamma(1/2)=\sqrt\pi$ and
$\Gamma(x+1)=x\Gamma(x)$ so $\Gamma(1+1/2)=\sqrt\pi/2$. Thus, with
$x=1/2$,
\[(\Gamma
\left(1+1/2\right))^{-k-\ell}
=\left(\frac{2}{\sqrt\pi}\right)^{k+\ell},
\]
so
\begin{eqnarray*}
a(k,\ell,1)=~~~~~~~~~~~~~~~~~~~~~~~~~~~~~~~~~~~~~~~~~~~~~~~~~~~~~~~~~~~~~~~~~~~~~~~~~~~~~~~~~~~~~~~~~~~~~~~
~~~~~~~~~~~~~~~~~~~~~~~\\
=\left[
\left(\frac{1}{\sqrt{2\pi}}\right)^{k+\ell-1}\cdot\left(\frac{1}{2}
\right)^{k\ell}\cdot(k+\ell)^{(k^2+\ell^2)/2}
\right]\cdot~~~~~~~~~~~~~~~~~~~~~~~~~~~~~~~~~~~~~~~~~~~~~~~~~~~~~~~~~~~~~~~~~~~~~~\\
\left(\frac{1}{k+\ell}
\right)^{\frac{1}{4}\cdot[k(k+1)+\ell(\ell+1)]}\cdot
\frac{1}{k!\cdot\ell!}\cdot
\frac{1}{\sqrt{2\pi}}\cdot
\left(\sqrt{2\pi} \right)^{k+\ell}
\cdot~~~~~~~~~~~~~~~~~~~~~~~~~~~~~~~~~~~~~~~~~~~~~~\\
\left(\frac{2}{\sqrt\pi}\right)^{k+\ell}\cdot
\prod_{i=1}^k\Gamma\left(1+i/2\right)\cdot
\prod_{j=1}^\ell\Gamma\left(1+j/2\right)= ~~~~~~~~~~~~~~~~~~\\
=\left(\frac{1}{\sqrt 2}
\right)^{2k\ell-2k-2l-1}\cdot\left(\frac{1}{\sqrt \pi}
\right)^{k+\ell-1}\cdot(k+\ell)^{\frac{1}{4}\cdot[k(k-1)+\ell(\ell-1)]}\cdot
\frac{1}{k!\cdot\ell!}
\cdot\frac{1}{\sqrt{2\pi}}\cdot~~~~~~~~~~~~~~~~~~~\\
\cdot\prod_{i=1}^k\Gamma\left(1+ i/2\right)\cdot
\prod_{j=1}^\ell\Gamma\left(1+j/2\right).~~~~~~~~~~~~~
\end{eqnarray*}
The factor $1/\sqrt{2\pi}$ cancels and we have
\begin{thm}\label{final2}
\[
S_{k,\ell}^{(1)}(n)\simeq
a(k,\ell,1)\cdot\left(\frac{1}{n}\right)^{\frac{1}{4}\cdot(k(k-1)+\ell(\ell-1))}\cdot
(k+\ell)^{n}
\]
where
\begin{eqnarray*}
a(k,\ell,1) =\left(\frac{1}{ 2}
\right)^{k\ell-k-\ell}\cdot\left(\frac{1}{\sqrt \pi}
\right)^{k+\ell}\cdot(k+\ell)^{\frac{1}{4}\cdot[k(k-1)+\ell(\ell-1)]}\cdot
\frac{1}{k!\cdot\ell!}\cdot \prod_{i=1}^k\Gamma\left(1+
i/2\right)\cdot
\prod_{j=1}^\ell\Gamma\left(1+j/2\right).~~~~~~~~~~~~~~~~~~~~~~~~
\end{eqnarray*}
\end{thm}

\subsubsection{A case with $2z=1$}

1. ~$k=\ell=1$. Then by Theorem~\ref{final2} $a(1,1,1)=1/2$ and $g(1,1,1)=0$, so
\[
{S^{(1)}_{1,1}}(n)\simeq {\frac{1}{2}}\cdot 2^n.
\]
On the other hand we know that
\[
{S^{(1)}_{1,1}}(n)=\sum_{j=0}^{n-1}{n-1\choose j}=2^{n-1},
\]
which verifies  Theorem~\ref{final2} (or Theorem~\ref{intro1})
in that case.

\subsubsection{Using "Mathematica"}\label{math1}
 For small $k$ and $\ell$ it is possible to write an explicit formula for, say, $S_{k,\ell}^{(1)}(n)$. By Theorem~\ref{intro1}
$S_{k,\ell}^{(1)}(n)\simeq A(k,\ell,n)$. Now form the ratio
$S_{k,\ell}^{(1)}(n)/ A(k,\ell,n)$. Using, say, "Mathematica", calculate that ratio
for increasing values of $n$, verifying that these values become closer and closer to 1 as $n$ increases. This indicates the validity  of Theorem~\ref{intro1}. We demonstrate this in the case $k=2$ and $\ell=1$.

\medskip
Here
$\frac{1}{4}\cdot(k(k-1)+\ell(\ell-1))=1/2$ and $a(2,1,1)=
\frac{1}{4}\cdot\sqrt{\frac{3}{\pi}} $. Thus by Theorem~\ref{intro1} (or~\ref{final2})
\[
S_{2,1}^{(1)}\simeq\frac{1}{4}\cdot\sqrt{\frac{3}{\pi}}\cdot
\frac{1}{\sqrt n}\cdot 3^n.
\]


 Next we  deduce a relatively simple formula for
 $S^{(1)}_{2,1}(n)$.
 By~\cite{regev3}
 \begin{eqnarray*}
S_{2,1}^{(1)}(n)=S(2,1;n)=~~~~~~~~~~~~~~~~~~~~~~~~~~~~~~~~~~~~~~~~~~~~~~~~~~~~~~~~~~~~~~~~~~~~~~~~~~~~~~~~~~~~~~~
~~~~~~\\=\frac{1}{4}\left(\sum_{r=0}^{n-1}{n-r\choose{\lfloor\frac{n-r}{2}\rfloor}}
{n\choose r}
+\sum_{k=1}^{\lfloor\frac{n}{2}\rfloor-1}\frac{n!}{k!\cdot
(k+1)!\cdot (n-2k-2)!\cdot (n-k-1)\cdot(n-k)}\right)+1.
\end{eqnarray*}
Also for $n\ge 2$ it can be proved by the WZ
method~\cite{doron1},~\cite{doron2} that
\[
2\sum_{j\ge 1}{n\choose j}{n-j\choose
j}=~~~~~~~~~~~~~~~~~~~~~~~~~~~~~~~~~~~~~~~~~~~~~~~~~~~~~~~~~~~~~~~~~~~~~~~~~~~~~~~~
\]
\begin{eqnarray}\label{motzkin.path.222}
=\sum_{r=0}^{n-1}{n-r\choose{\lfloor\frac{n-r}{2}\rfloor}}
{n\choose r}
+\sum_{k=1}^{\lfloor\frac{n}{2}\rfloor-1}\frac{n!}{k!\cdot
(k+1)!\cdot (n-2k-2)!\cdot (n-k-1)\cdot(n-k)}
\end{eqnarray}
(for an elementary proof of Equation~\eqref{motzkin.path.222} (due
to I. Gessel), see~\cite{regev4}). Hence
\[
S_{2,1}^{(1)}\simeq \frac{1}{2}\cdot\sum_{j\ge 1}{n\choose
j}{n-j\choose j}
\]
If indeed
\[
S_{2,1}^{(1)}\simeq\frac{1}{4}\cdot\sqrt{\frac{3}{\pi}}\cdot
\frac{1}{\sqrt n}\cdot 3^n
\]
then
\begin{eqnarray}\label{correct}
\frac{2\cdot\sqrt{n}}{3^n}\cdot\sum_{j\ge 1}{n\choose
j}{n-j\choose j}\simeq\sqrt{\frac{3}{\pi}}=0.977205.
\end{eqnarray}
Indeed, "Matematica" gives the following values $lhs(n)$ for the
left hand side of~\eqref{correct}:

\medskip
\[lhs(10)=0.958821, ~lhs(100)=0.975373, ~lhs(1000)=0.977022, \]

\[~lhs(2000)=0.977113, ~lhs(3000)= 0.977144\quad\mbox{ etc.}\]
This, in a sense, verifies Theorem~\ref{intro1} in this case.

\subsubsection{A special case with $z=1$}

Let $k=\ell=z=1$. By Theorem~\ref{final1}
\begin{eqnarray}\label{sof1}
S_{1,1}^{(2)}(n)\simeq a\cdot \left(\frac{1}{n}\right)^g\cdot
2^{2n}\quad\mbox{ where}\quad a=\frac{1}{4\sqrt{\pi}}\quad\mbox{
and}\quad g=\frac{1}{2}.
\end{eqnarray}
Now
\[
S_{1,1}^{(2)}(n)=\sum _{j\ge 0}{n-1\choose j}^2={2(n-1)\choose
n-1}
\]
and by Stirling's formula
\[
{2(n-1)\choose n-1}\simeq \frac{1}{\sqrt\pi}\cdot\frac{1}{\sqrt
n}\cdot 4^{n-1},
\]
agreeing with~\eqref{sof1}, hence again, Theorem~\ref{intro1} is verified in this case.

A. Berele, Department of Mathematics, DePaul University, Chicago,
Ill 60614

email: aberele~at~condor.depaul.edu

\medskip
A. Regev, Department of Mathematics, The Weizmann Institute of
Science, Rehovot 76100,  Israel

e-mail: amitai.regev~at~weizmann.ac.il


\begin{thebibliography}{99}

\bibitem{berele} A. Berele and A. Regev, Hook Young diagrams with
applications to combinatorics and to representation theory of Lie
superalgebras, Advances in Math. Vol. 64 No. 2 (1987) 118-175

\bibitem{forrester} P. J. Forrester and S. O. Warnaar, The importance of the Selberg
integral. Bull. Amer. Math. Soc. (N.S.) 45 (2008), no. 4,
489--534.

\bibitem{garsia} A. M. Garsia and N. Wallach, The non-degeneracy of the bilinear
form of $m$-Quasi-Invariants, Advances in Applied Math, vol. {\bf
37}, (2006) 309--360.


\bibitem{doron1} M. Petkovesk, H.S. Wilf and D Zeilberger, A=B, AK
Peters Ltd. (1996)


\bibitem{regev} A. Regev, Asymptotic values for degrees associated
with strips of Young diagrams, Advances in Math. {\bf 41}, No 2
(1981) 115-136

 \bibitem{regev3} A. Regev, Probabilities in the $(k,\ell)$ hooks,
 Israel J. Math. {\bf 169} (2009) 61-88

 \bibitem{regev2} A. Regev,  Humps for Dyck and for Motzkin paths,
 Preprint (2010).

\bibitem{regev4} A. Regev, Identities for the number of standard Young tableaux in some
$(k,\ell)$ hooks, arXiv


\bibitem{selberg} A. Selberg, Bemerkninger om et multipelt
integral, {\it Norsk Mat. Tidskr.} {\bf 26} (1944), 71-78.

 \bibitem{doron2} D. Zeilberger, The method of creative telescoping  J. Symbolic Computation
 11, 195-204 (1991).


\end{thebibliography}
\end{document}